\documentclass{article}
\usepackage{row}
\usepackage{graphicx}
\usepackage[utf8]{inputenc}
\begin{document}
\title{The Geometry of Manifolds and the Perception of Space}
\author{Raymond O. Wells, Jr.
\footnote{Jacobs University Bremen; Rice University; University of Colorado at Boulder; row@raw.com}}
\maketitle

\section{Introduction}
In our contemporary world, we see ourselves immersed in a vast universe filled with galaxies, black holes, dark matter and other aspects of our cosmological surroundings. An important aspect to all of this is that we imagine ourselves to be in a  variably curved four-dimensional space-time, as first put forward by Einstein in his theory of general relativity in 1916 \cite{einstein1916}.   In February of 2016, a century after Einstein's paper appeared, it was announced \cite{abbott2016} by a large team of astronomers that gravitational waves were first detected by amazingly accurate sensors. This result produced headlines around the world, as did the confirmation of Einstein's theory by Dyson, Eddington, and Davidson in their celebrated eclipse paper \cite{dyson-eddington-davidson1920} of 1920 showing that light was bent by the gravitational attraction of the sun, as was predicted by Einstein. The measurement of the bending of a specific star was \(1.65''\) which compared favorably with Einstein's prediction of \(1.75''\).

This four-dimensional space-time is a four-dimensional manifold with a Lor\-ent\-zian metric which satisfies the Einstein field equations.  The question we want to discuss in this essay is: what are abstract manifolds (four-dimensional or otherwise), and how did this mathematical concept and its associated geometric ideas arise, so that it could be used for modeling our universe as Einstein and so many others after him have done?

\section{The Origins of the Concept of a Manifold}
\label{sec:origins-manifold}
The notion of a manifold is a relatively recent one, but the theory of curves and surfaces in Euclidean three-space \(\BR^3\) originated in the Greek mathematical culture. For instance, the book on the study of conic sections  by Apollonius  \cite{apollonius1896} described mathematically the way we understand conic sections today.  The intersection of a plane with a cone in \(\BR^3\) generated these curves, and Apollonius showed moreover that any plane intersecting  a {\it skew cone} gave one of the three classical conic sections (ellipse, hyperbola, parabola), excepting the degenerate cases of a point or intersecting straight lines. This was a difficult and important theorem at the time. The Greek geometers studied intersections of other surfaces as well, generating additional curves useful for solving problems and they also introduced some of the first curves defined by transcendental functions (although that terminology was not used at the time), e.g., the quadratrix, which was used to solve the problem of squaring the circle (which they couldn't solve with straight edge and compass and which was shown many centuries later not to be possible).  See for instance Kline's very fine book on the history of mathematics for a discussion of these issues \cite{kline1972}.

Approximately 1000 years after the major works by the Greek geometers, Descartes published in 1637 \cite{descartes} a revolutionary book which contained a fundamentally new way to look at geometry, namely as the solutions of algebraic equations.  In particular the solutions of equations of second degree in two variables described precisely the conic sections that Apollonius had so carefully treated.  This new bridge between algebra and geometry became known during the 18th and 19th centuries as {\it analytic geometry} to distinguish it from {\it synthetic geometry}, which was the treatment of geometry as in the book of Euclid \cite{euclid}, the methods of which were standard in classrooms and in research treatises for approximately two millenia. Towards the end of the 19th century and up until today {\it algebraic geometry} refers to the relationship between algebra and geometry, as initiated by Descartes, as other forms of geometry had arisen to take their place in modern mathematics, e.g., topology, differential geometry, complex manifolds and spaces, and many other types of geometries.

Newton and Leibniz made their major discoveries concerning differential and integral calculus in the latter half of the 17th century (Newton's version was only published later; see \cite{newton1736}, a translation into English by John Colson from a Latin manuscript from 1671). Leibniz published his major work on calculus in 1684 \cite{leibniz1684}.  As is well known, there was more than a century of controversy about priority issues in this discovery of calculus  between the British and Continental scholars (see, e.g., \cite{kline1972}). In any event, the 18th century saw the growth of analysis as a major force in mathematics (differential equations, both ordinary and partial, calculus of variations, etc.). In addition the variety of curves and surfaces in \(\BR^3\) that could be represented by many new transcendental functions expanded greatly the families of curves and surfaces in \(\BR^3\) beyond those describable by solutions of algebraic equations. 

Moreover, a major development that grew up at that time, and which concerns us in this paper, was the growing interaction betwen analysis and geometry. An important first step was the analytic description of the curvature of a curve in the plane at a given point by Newton. This was published in Newton's 1736 monograph  \cite{newton1736} mentioned earlier.  In the text of this reference (in a section entitled ``To find the quantity of curvature in any curve'') one finds the well-known formula for the curvature of a curve defined as the graph of a function \(y=f(x)\)
\be
\label{eqn:plane-curvature}
K_P=\pm \frac{f''(x)}{[1+(f'(x)^2]^\frac{3}{2}},
\ee
which one learns early on in the study of calculus.  Curvature had been studied earlier quite extensively by Huygens \cite{huygens1673}, and he computed the curvature of many explicit curves (including the conic sections, cycloids, etc.). He did not  use calculus,{\it per se}, but he used limiting processes of geometric approximations that are indeed the essence of calculus. Even earlier Apollonius had been able to compute the curvature of a conic section at a specific point (see \cite{apollonius1896} and Heath's discussion of this point).

The next development concerning analysis and geometry was the study of space curves in 
\(\BR^3\) with its notion of curvature and torsion as we understand it today.  This work on space curves was initiated by a very young (16 years old) Clairaut \cite{clairaut1731} in 1731.  Over the course of the next century there were numerous contributions to this subject by Euler, Cauchy, and many others, culminating in the definitive papers of Frenet \cite{frenet1852} and Serret \cite{serret1851} in the mid-19th century giving us the formulas for space curves in \(\BR^3\) that we learn in textbooks today.  Space curves were called ``curves of double curvature" (French: {\it courbes \`{a} double courbure}) throughout the 18th and 19th centuries.

A major aspect of studying curves in two or three dimensions is the notion of the measurement of arc length, which has been studied since the time of Archimedes, where he gave the first substantive approximate value of \(\pi\) by approximating the arc length of an arc of a circle (see, e.g., \cite{resnikoff-wells1984}).  Using calculus one was able to formalize the length of a segment of a curve by the well-known formula
\[
\int_\G ds = \int^b_a \sqrt{(x'(t))^2+(y'(t))^2}dt,
\]
if \(\G\) is a curve in \(\BR^2\) parametrized by \((x(t),y(t))\) for \(a\le t \le b\). The expression \(ds^2= dx^2 + dy^2\) in \(\BR^2\) or \(ds^2 =dx^2 +dy^2 +dz^2\) in \(\BR^3\) became known as the {\it line element} or infinitesimal measurement of arc length for any parametrized curve.  What this meant is that a given curve could be approximated by a set of secants, and measuring the lengths of the secants by the measurement of the length of a segment of a straight line in \(\BR^2\) or \(\BR^3\) (using Pythagoras!), and taking a limit, gave the required arc length. The line elements \(ds^2= dx^2 +dy^2\) and \(d^2 = dx^2 +dy^2 +dz^2\) expressed {\it both} the use of measuring straight-line distance in the ambient Euclidean space as well as the limiting process of calculus.  We will see later how generating such a line element to be {\it independent} of an ambient Euclidean space became one of the great creations of the 19th century.

The concept of curvature of a curve in \(\BR^3\) was well understood at the end of the 18th century, and the later work of Cauchy, Serret and Frenet completed this set of investigations begun by the young Clairaut a century earlier.  The problem arose:  how can one define the {\it curvature of a surface} defined either locally or globally in \(\BR^3\)? 

Euler recognized the difficulty of this problem, as we see from this quote of Euler's from his paper on curvature of surfaces from 1767,   ``Recherches sur la courbure des surfaces" \cite{euler1767}:
\begin{quote}
In order to know the curvature of a curve, the determination of the radius of the osculating circle furnishes us the best measure, where for each point of the curve, we find a circle whose curvature is precisely the same. However, when one looks for the curvature of a surface, the question is very equivocal and not at all susceptible to an absolute response, as in the case above. There are only spherical surfaces where one would be able to measure the curvature, assuming the curvature of the sphere is the curvature of its great circles, and whose radius could be considered the appropriate measure. But for other surfaces, one doesn't know even how to compare a surface with a sphere, as when one can always compare the curvature of a curve with that of a circle. The reason is evident, since at each point of a surface there are an infinite number of different curvatures. One has to only consider a cylinder, where along the directions parallel to the axis, there is no curvature, whereas in the directions perpendicular to the axis, which are circles, the curvatures are all the same, and all other obliques sections to the axis give a particular curvature. It's the same for all other surfaces, where it can happen that in one direction the curvature is convex, and in another it is concave, as in those resembling a saddle.
\end{quote}

In the quote above, we see that Euler recognized the difficulties in defining curvature for a surfacs at any given point.  He does not resolve this issue in this paper, but he makes extensive calculations and several major contributions to the subject. In particular, he computes the curvatures at a point of the surface of the curves which are the intersections of planes passing though a normal vector at the point with the surface and finds what are now called the {\it principal curvatures} and the {\it principal directions} of the surface at the given point.

\section{Gauss and Intrinsic Differential Geometry}
\label{sec:gauss}

A major development at the beginning of the 19th century was the publication in 1828 by Gauss  of his historic landmark paper concerning differential geometry of surfaces, entitled {\it Disquisitiones circa superficies curvas} \cite{gauss1828}.  The year before, he published a very readable announcement and summary of his major results in \cite{gauss1827}, and we shall quote from this announcement paper somewhat later, letting Gauss tell us in his own words what he thinks the significance of his discoveries is.  For the moment, we will simply say that this paper laid the foundation for doing intrinsic differential geometry on a surface and was an important step in the creation of a theory of abstract manifolds which was developed a century later.

Gauss defined {\it curvature} \(\kappa\) of a surface \(S\subset\BR^3\) to be the derivative of what is now called the Gauss mapping, which maps a point on the surface to its unit normal vector, which is a point on the unit sphere in \(\BR^3\). He showed first that the curvature (now called the {\it Gaussian curvature}) is the product of the Eulerian prinipal curvatures at that point. This definition of curvature (either via the Gauss mapping or the product of the principal curvatures) depends explicitly on the embedding of the surface in \(\BR^3\).  They both use the notion of the normal vector at the given point, which depends very much on the embedding.

Suppose that the surface is described in terms of local coordinates \((u,v)\), i.e.,
\be
\label{eqn:surface-coordinates}
x(u,v), y(u,v), z(u,v)
\ee
are three smooth functions that represent the surface parametrically, then the line element 
\[
ds^2=dx^2+dy^2+dz^2
\]
in \(\BR^3\) induces a metric on the surface \(S\) of the form
\be
\label{eqn:gauss-metric}
ds^2=Edu^2+2F dudv+ Gdv^2,
\ee
where \(E,F,\) and \(G\) are functions of \((u,v)\) obtained by differentiating the para\-met\-rizing functions \(x(u,v)\), \(y(u,v)\), and \(z(u,v)\).  One can use the line element \(ds^2\) on \(S\) to measure the length of a curve on the surface, to measure the angle formed by two intersecting curves on \(S\), and to compute other geometric quantities.

The great achievement of Gauss in this paper was to show (his {\it Theorema Egregium}) that the Gaussian curvature could be defined in terms of the first and second derivatives of the coefficients \(E, F\), and \(G\) of the line element in (\ref{eqn:gauss-metric}), and was thus independent of the embedding. A major motivation for Gauss's work on curvature in this paper were his parallel experimental geodesic measurements, which were trying to estimate the curvature of the earth near Göttingen.

Gauss understood full well the significance of his work and the fact that this was the beginning of the study of a new type of geometry (which later generations have called {\it intrinisic differential geometry}).  We quote here from his announcement of his results published some months earlier from pp. 344--345 of \cite{gauss1827}:
\begin{quote}
\label{gauss-quote}
These theorems lead us to consider the theory of curved surfaces from a new point of view, whereby the investigations open to a quite new undeveloped field.  If one doesn't consider the surfaces as boundaries of domains, but as domains with one vanishing dimension, and at the same time as bendable but not as stretchable, then one understands that one needs to differentiate between two different types of relations, namely, those which assume the surface has a particular form in space and those that are independent of the different forms a surface might take.  It is this latter type that we are talking about here.  From what was remarked earlier, the curvature belongs to this type of concept; and, moreover, figures constructed on the surface, their angles, their surface area, their total curvature as well as the connecting of points by curves of shortest length, and similar concepts, all belong to this class.
\end{quote}
The results described in this short announcement (and the details in the much longer Latin paper on the subject) formed the basis of most of what became modern differential geometry.  

An important point that we should make here is that Gauss did significant experimental work on measuring the curvature of the earth in the area around G\"{o}ttingen, where he spent his whole scientific career.  This involved measurements over hundreds of miles, and involved communicating between signal towers from one point to another.  He developed his theory of differential geometry as he was conducting the experiments, and at the end of his announcement \cite{gauss1827}, he summarizes one of his experiments to say that a geodesic triangle with the longest side 12 miles long has the sum of its three angles measuring \(2''\) greater that 180\(^\circ\), indicating in a precise manner the positive curvature of the earth near that triangle. It is interesting to compare these experiments with those concerning the bending of light one century later and the measurement of gravitional waves two centuries later.

\section{Riemann's Higher-Dimensional Geometry}
In mathematics, we sometimes see striking examples of brilliant contributions or completely new ideas that change the ways mathematics develops in a significant fashion.   A prime example of this is the work of Descartes \cite{descartes}, which completely changed how mathematicians looked at geometric problems.  But it is rare that a single mathematician makes as many singular advances in his lifetime as did Riemann in the middle of the 19th century. In this section, we will discuss in some detail his fundamental creation of the theory of higher-dimensional manifolds and the additional creation of what is now called Riemannian or simply differential geometry.  However, it is worth noting that he only published nine papers in his short lifetime (he lived to be only 40 years old), and several other important works, including those that concern us in this section, were published posthumously from the writings he left behind.  His collected works (including in particular these posthumously published papers) were edited and published in 1876 and are still in print today \cite{riemann1876}.

Looking through the titles, one is struck by the wide diversity as well as the originality.  Let us give a few examples here.  In Paper I%
\footnote{We use the enumeration from the Table of Contents of his collected works \cite{riemann1876}, which include the posthumously published papers as well.}
 (his dissertation), he formulated and proved the Riemann mapping theorem and dramatically moved the theory of functions of one complex variable in new directions.  In Paper VI, in order to study Abelian functions, he formulated what became known as Riemann surfaces and this led to the general theory of complex manifolds in the 20th century.  In Paper VII, he proved the Prime Number Theorem and formulated the Riemann Hypothesis, which is surely the outstanding mathematical problem in the world today. In Paper XII, he formulated the first rigorous definition of a definite integral (the Riemann integral) and applied it to trigonometric series, setting the stage for Lebesgue and others in the early 20th century to develop many consequences of the powerful theory of Fourier analysis.  In Papers XIII and XXII, he formulated the theory of higher-dimensional manifolds, including the important concepts of Riemannian metric, normal coordinates, and the Riemann curvature tensor, which we will visit very soon in the paragraphs below.  Paper XVI contains correspondence with Enrico Betti leading to the first higher-dimensional topological invariants beyond those already known for two-dimensional manifolds (that had been developed by Riemann in Papers I and VI).

His paper \cite{riemann1854} (Paper XIII above) is a posthumously published version of a public lecture Riemann gave as his { \it Habilitationsvortrag} in 1854. This was part of the process for obtaining his {\it Habilitation}, a German advanced degree beyond the doctorate necessary to qualify for a professorship in Germany at the time (such requirements still are in place at most German universities today as well as in other European countries, e.g., France and Russia; it is similar to the research requirements in the US to be qualified for tenure). This paper, being a public lecture, has very few formulas, is at times quite philosophical and is amazing in its depth of vision and clarity.  On the other hand, it is quite a difficult paper to understand in detail, as we shall see.

Before this paper was written, manifolds were all one- or two-dimensional curves and surfaces in \(\BR^3\), including their extension to points at infinity (which was developed as a part of projective geometry in the first half of the 19th century).  In fact, some mathematicians who had to study systems parametrized by more than three variables declined to call the parametrization space a manifold or give such a parametrization a geometric significance.  In addition, these one- and two-dimensional manifolds always had a differential-geometric structure which was induced by the ambient Euclidean space (this was true for Gauss, as well).

In Riemann's paper \cite{riemann1854}, he discusses the distinction between discrete and continuous manifolds, where one can make comparisons of quantities by either counting or by measurement, and gives a hint, on p. 256, of the concepts of set theory, which was only developed later in a single-handed effort by Cantor. Riemann begins his discussion of manifolds by moving a one-dimensional manifold, which he intuitively describes, in a transverse direction (moving in some type of undescribed ambient ``space"), and inductively, generating an \(n\)-dimensional manifold by moving an \((n-1)\)-manifold transversally in the same manner.  Conversely, he discusses having a nonconstant function on an \(n\)-dimensional manifold, and the set of points where the function is constant is (generically) a lower-dimensional manifold; and by varying the constant, one obtains a one-dimensional family of \((n-1)\)-manifolds (similar to his construction above).%
\footnote{He alludes to some manifolds that cannot be described by a finite number of parameters; for instance, the manifold of all functions on a given domain, or all deformations of a spatial figure.  Infinite-dimensional manifolds, such as these, were studied in great detail a century later.}

Riemann formulates local coordinate systems \((x^1, x^2,...,x^n)\) on a manifold of \(n\) dimensions near some given point, taken here to be the origin. He formulates a curve in the manifold as being simply \(n\) functions \((x^1(t), x^2(t),...,x^n(t))\) of a single variable \(t\). The concepts of set theory and topological space were developed only later in the 19th century, and so the global nature of manifolds is not really touched on by Riemann (except in his later work on Riemann surfaces and his correspondence with Betti, mentioned above). It seems clear on reading his paper that he thought of \(n\)-dimensional manifolds as being extended beyond Euclidean space in some manner, but the language for this was not yet available.

At the beginning of this paper Riemann acknowledges the difficulty he faces in formulating his new results.  Here is a quote from the second page of his paper (p. 255):
\begin{quote}In that my first task is to try to develop the concept of a multiply spread-out quantity [he uses the word manifold later], I believe even more in  being allowed an indulgent evaluation, as in such works of a philosophical nature, where the difficulties are more in the concepts than in the construction, wherein I have little experience, and except for the paper by Mr. Privy Councilor Gauss in his second commentary on biquadratic residues in the G\"{o}ttingen Gelehrte Anzeige [1831]  and in his J{u}bil\"{a}umsschrift and some investigations by Hebart, I have no precedents I could use.
\end{quote}
The paper of Gauss that he cites here \cite{gauss1831} refers to Gauss's dealing with the philosophical issue of understanding the complex number plane after some thirty years of experience with its development.  We will mention this paper again somewhat later in this paper.  Hebart  was a philosopher whose metaphysical investigations influenced Riemann's thinking.  Riemann was very aware of the speculative nature of his theory, and he used this philosophical point of view, as the technical language he needed (set theory and topological spaces) was not yet available.  This was very similar to Gauss's struggle with the complex plane, as we shall see later.

As mentioned earlier, measurement of the length of curves goes back to the Archimedian study of the length of a circle.  The basic idea there and up to the work of Gauss was to approximate a given curve by a collection of straight line segments and take a limit.  The {\it length} of each straight line segment was determined by the Euclidean ambient space, and the formula, using calculus for the limiting process, became, in the plane for instance,
\[
\int_\G ds = \int^b_a\sqrt{(x'(t))^2 + (y'(t))^2}dt,
\]
where \(ds^2=dx^2+dy^2\) is the line element of arc length in \(\BR^2\).  As we saw in Section \ref{sec:gauss}, Gauss formulated in \cite{gauss1828} on a two-dimensional manifold with coordinates \((p,q)\) the line element
\be
\label{eqn:line-element}
ds^2=Edp^2 + 2Fdpdq + Gdq^2
\ee
where \(E,F,\) and \(G\) are induced from the ambient space.  He didn't consider any examples of such a line element (\ref{eqn:line-element}) that weren't induced from an ambient Euclidean space, but his remarks (see the quote above in Section \ref{sec:gauss}) clearly indicate that this could be a ripe area for study, and this could well include allowing coeffients of the line element  (\ref{eqn:line-element}) being more general than induced from an ambient space.

Since Riemann formulated an abstract \(n\)-dimensional manifold (with a local coordinate system) with no ambient space, and since he wanted to be able to measure the length of a curve on his manifold, he formulated, or rather postulated, an independent measuring system which mimics Gauss's formula (\ref{eqn:line-element}).  Namely, he prescribes for a given coordinate system a metric (line element) of the form 
\be
\label{eqn:riemann-metric}
ds^2= \sum^n_{i,j=1} g_{ij}(x)dx^i dx^j,
\ee
where  \(g_{ij}(x) \) is, for each \(x\), a symmetric positive definite matrix, and he postulates by the usual change of variables formula
\[
ds^2= \sum^n_{i,j=1}\tilde{g}_{ij}(\tilde{x}) d \tilde{x}^i d \tilde{x}^j,
\]
where \(\tilde{g}_{ij}(\tilde{x})\) is the transformed positive definite matrix in the new coordinate system \((\tilde{x}_1,...,\tilde{x}^n
)\). 

Using the line element (\ref{eqn:line-element}), the length of a curve is defined by 
\be
\label{eqn:riemannian-metric}
l(\G):= \int^b_a\sqrt{\sum^n_{i,j=1} g_{ij}(x(t))\frac{dx^i}{dt}(t)\frac{dx^j}{dt}(t)}dt.
\ee
The line element (\ref{eqn:riemann-metric}) is what is called a {\it Riemannian metric} today, and the 2-form \(ds^2\) is considered as a positive definite bilinear form giving an inner product on the tangent plane \(T_p (M)\) for \(p\) a point on the manifold \(M\). This has become the basis for almost all of modern differential geometry (with the extension to Lorentzian-type spaces where \(g_{ij}(x)\) is not positive definite, {\it \`{a} la} Minkoswki space). Riemann merely says
on page 260 of his paper (no notation here at all),
\begin{quote}
I restrict myself therefore to manifolds where the line element is expressed by the square root of a differential expression of second degree.
\end{quote}
Earlier he had remarked that a line element should be homogeneous of degree one and one could also consider the fourth root of a differential exression of fourth degree, for instance.  Hence his restriction in the quote above.

The next step in Riemann's paper is his formulation of curvature.  This occurs on a single page (p. 261 of \cite{riemann1854}), 
It is extremely dense and not at all easy to understand.  Over time, however, it became understood by a several generations of mathematicians, and this became the basis for the work of Einstein.

We summarize briefly what Riemann described on this page. He formulated the notion of geodesic coordinates, i.e., coordinates \((x_1,\ldots,s_n)\), where the coordinates are geodesics in the metric \(ds^2\) from (\ref{eqn:riemannian-metric}), and he expanded the functions \(g_{ij}\) in a Taylor series at the origin and considered the second-order terms as a biquadratic form 
\be
\label{eqn:riemann-curvature1}
Q(x,dx)=\sum_{ijkl}c_{ijkl} x^kx^l dx^i dx^j,
\ee
and the coefficients \(c_{ijkl}\) had certain important symmetry properties and  effectively defined the {\it Riemannian curvature tensor} written classically and still today as a tensor of the form
\[
R^\r_{\s\m\n}.
\]

How does one define such a curvature tensor for \(n\)-dimensional manifolds with a Riemannian metric in a general coordinate system? In a quite short  paper written in Latin for a particular mathematical prize in Paris in 1861, (Paper No. XXII) Riemann provides the first glimpse of the general Riemann curvature tensor.  The purpose of Riemann's paper was to answer a question in the Paris competition dealing with the flow of heat in a homogeneous solid body.

Riemann's ideas in these two papers were developed and expanded considerably in the following decades in the work of Christoffel, Levi-Cevita, Ricci, Beltrami, and many others. The main point of our discussion has been that Riemann created on these few pages the basic idea of an \(n\)-dimensional manifold not considered as a subset of Euclidean space {\it and} of the independent concept of a Riemannian metric and the Riemannian curvature tensor.  What is missing at this point in time is the notion of a topological space on the  basis of which one can formulate the contemporary concept of an abstract manifold or an abtract Riemannian manifold.

\section{Hermann Weyl: Prelude to the 20th Century}
  Other major developments in geometry that led to the study of abstract manifolds in the 20th century, and which we won't consider in any detail in this paper, are the creations of complex geometry,  transformations groups, set theory and  topology, among others. Complex numbers were called {\it imaginary numbers} (among other appellations) in the centuries  preceding the 19th century.  This meant, in the eyes of the beholders that these were not concrete real (or realistic) numbers (to use a pun!), but were simply imaginary artifacts that had no real meaning, but were useful in the way they arose as the would-be solutions of algebraic equations. Gauss, among others, used these numbers extensively in his career, for instance in his ground-breaking work on number theory early in his life, {\it Disquisitiones Arithmeticae} \cite{gauss1801}. Only many decades later did he make the case for a definitive geometric interpretation of what we call complex numbers today in a brief paper \cite{gauss1831}, which was a commentary on some of his earlier work on number theory and which Riemann cited earlier in this paper as a philosophical work which was a guide for him. Gauss called these numbers (for the first time) {\it complex numbers}, with their real and imaginary parts representing coordinates in a two-dimensional plane --- the {\it complex plane}. Gauss pleads with his readers to consider complex numbers and the complex plane to be considered as a normal part of mathematics, not as something ``imaginary" or ``unreal."  It's very enlightening to read this quite readable short paper.  Later in the century, Klein and others started to consider complex solutions of homogeneous algebraic equations and thus began the study of complex algebraic manifolds and varieties.  Klein along with Lie fostered the study of transformation groups, and the notion of manifolds as quotients of such groups became an important development as well.  See, for instance, the very informative book by Klein \cite{klein1928} (first published in 1928, but based on lectures of Klein towards the end of the 19th century), which described most of the developments in geometry in the 19th century in a succinct fashion.

The creation of set theory in the latter half of the 19th century was a singular effort of Georg Cantor over several decades of work.  The first article in this direction was published in 1874 \cite{cantor1874}. He then published six major papers in the Mathematische Annalen between 1879 and 1884,
which established the basic tenets of set theory, and laid the foundation for work in multiple directions for the next century and beyond, in particular in logic and foundations of mathematics, as well as point set topology, to mention a particular part of geometry related to our thesis in this paper (see his collected works \cite{cantor1932} for these papers).

After  Cantor's groundbreaking discovery (some aspects of which were very controversial for quite some time) the notion of topological space evolved to its contemporary state as a set satisfying certain axioms concerning either neighborhoods or open sets.  This became the development of point set topology, which has its own very interesting history, and we mention only the fundamental contributions of Fréchet in 1906 \cite{frechet1906} (metric spaces, Fréchet spaces, in particular in the infinite-dimensional setting) and Hausdorff in 1914 \cite{hausdorff1914} (which, among other things, defined specifically Hausdorff spaces, one of many specialized types of topological spaces used often in geometry today). In 1895 Poincar\'{e} launched the theory of topological manifolds with a cornerstone paper called {\it Analysis Situs} \cite{poincare1895}, which was followed up by five supplements to this work over the next 17 years (see Vol. 6 of his collected works \cite{poincare1953} as well as a very nice translation of all of these topology papers of Poincar\'{e} by John Stillwell \cite{poincare2010}). This paper and its supplements became the foundation of what is now called {\it algebraic topology} (following the pioneering work of Riemann \cite{riemann1857} on Riemann surfaces and Betti on homological invariants of higher-dimensional manifolds \cite{betti1871}).

The final step in our journey  is the book by Hermann Weyl in 1913 \cite{weyl1913} entitled {\it Die Idee der Riemannschen Fl\"{a}che},%
\footnote{{\it The Concept of a Riemann Surface}.}
in which he gives the first very specific definition of an abstract manifold with more structure than a topological manifold (which Poincar\'{e} and others had already investigated quite thoroughly).  His motivation was to give a better understanding of Riemann surfaces which transcended Riemann's original description as a multisheeted covering of the complex plane with branch points of various kinds. Specifically, he considered a topological manifold of two dimensions (locally homeomorphic to an open set in \(\BR^2\) ) which had a finite or countable triangulation and which had the additional property that there were coordinate charts mapping to an open ball in the complex plane \(\BC\) whose transition functions on overlapping coordinate charts were holomorphic. He showed how all of the previous work on Riemann surfaces fits into this new picture, and he proved a fundamental existence theorem of global holomorphic or meromorphic functions on such surfaces, utilizing the Dirichlet principle, that had been first used by Riemann in his conformal mapping theorem from his dissertation \cite{riemann1851}.  This definition of a topological manifold with such additional structure became the model for all the various kinds of manifolds studied in the following century up to the current time. This included, for instance, differentiable (\(C^\infty\)) manifolds, complex manifolds of arbitrary dimension, Riemannian manifolds, symplectic manifolds, and real-analytic manifolds, among many others. Weyl knew that this was new territory, and like Gauss with intrinsic differential geometry, and Riemann with \(n\)-dimensional manifolds, he carefully explained to his readers that he was introducing a new way of thinking. Here, in Weyl's own words, are what he thought about this (\cite{weyl1913}, p. V):
\begin{quote}
Such a rigorous presentation, which, namely by the establishing of the fundamental concepts and theorems in function theory and using theorems of the analysis situs which do not just depend on intuitive plausibility, but have set-theoretic exact proofs, does not exist. The scientific work that remains to be done   in this regard may perhaps not be particularly highly valued.  But, nevertheless, I believe I can maintain that I have tried in a serious and conscientious manner to find the simplest and most appropriate methods that lead to the  asserted goal; and at many points, I have had to proceed in a different manner than that which has become traditional in the literature since the appearance of C. Neumann's classical book ``Lectures on Riemann's theory of Abelian Integrals" (1865).
\end{quote}

As we mentioned in the Introduction, Einstein's paper utilizing a Riemann\-ian structure with a Lorentzian metric appeared in 1916.  Almost immediately thereafter in the summer of 1917, Hermann Weyl gave lectures in Zurich on Einstein's theory, which appeared in a book published by Springer soon thereafter. The book was entitled, in German, {\it Raum, Zeit, Materie}, and its third edition appeared in 1919 \cite{weyl1919} (an English edition, {\it Space, Time, Matter}, became available in 1950, Dover Publications, New York, and is still available today). Weyl gave more mathematical background for Einstein's theory and had a major influence on  the propagation of Einstein's ideas.  Weyl considers arbitrary \(n\)-dimensional Riemannian manifolds and the corresponding theory of tensors on such a manifold.  The global structure of such a manifold is not discussed; he concentrates primarily on the local theory for varying coordinate systems near a specific point.

If we jump almost a century later, we see that modern string theory considers the four-dimensional space-time of Einstein, which we can denote by \(M^4\), as a basis for a 10-dimensional string-theoretic manifold. This is a four-dimensional space-time \(M^4\) with an additional compact three-dimensional complex manifold \(X^3\) (a Calabi--Yau manifold) attached to \(M^4\) at each point, and which is extremely small relative to our usual perception of space around us.  See, for example, the very interesting exposition of Shing-Tung Yau and Steve Nadis \cite{yau-nadis2011} for a discussion of this and with many references to the string-theoretic literature.  The notion of an abstract manifold, as first formulated by Hermann Weyl in 1913, is essential for this theory. It is also interesting to note that Weyl's book dealt with the important example of Riemann surfaces, which were one-dimensional complex manifolds (with local holomorphic coordinate systems) with specific global proerties (for instance, connectivity, which Riemann had initiated). This is, in fact, part of the nature of Calabi--Yau manifolds, which are three-dimensional generalizations of the one complex-dimensional Riemann surfaces of Hermann Weyl, and which have specific global properties, which distinguishes them from from other three-dimensional complex manifolds.

\bibliography{references}
\bibliographystyle{plain}%

\end{document}